\documentclass[a4paper]{article}

\usepackage{amsmath,amssymb}
\usepackage{cite}
\usepackage{amsfonts,mathtools}
\usepackage{graphicx,subfigure}
\usepackage[mathscr]{eucal}
\usepackage{epsfig,color,pifont,enumerate,rotating}
\usepackage{geometry}
\DeclareMathAlphabet{\mathpzc}{OT1}{pzc}{m}{it}

\newcommand{\ts}{\mathfrak{T}}

\renewcommand{\epsilon}{\varepsilon}

\newcommand{\br}{\mathbb{R}}
 \newcommand{\bldr}{\boldsymbol{r}}
 
 \newcommand{\bldp}{\boldsymbol{p}}

  \newcommand{\bldb}{\boldsymbol{s}}

 \newcommand{\bcs}{\boldsymbol{c}}
  
\newcommand{\e}{\mathpzc{E}}

\newtheorem{define}{Definition}[section]

\newtheorem{assumption}{Assumption}[section]
\newtheorem{lemma}{Lemma}[section]
\newtheorem{corr}{Corollary}[section]
\newtheorem{proposition}{Proposition}[section]

\newcommand{\spr}[2]{\left\langle #1; #2 \right\rangle}

\newcommand{\ff}[3]{\mathbf{II}_{#1}\!\left[ #2; #3 \right]}

\newcommand{\ov}[1]{\overline{#1}}

\newcommand{\ve}{\varepsilon}

\newcommand{\pf}{\par {\bf Proof:}\;}
\newcommand{\epf}{$\qquad \Box$}
\newcommand{\dist}[2]{\text{\bf d}_{#2} \left[ #1\right]}

\newcommand{\bpi}{\boldsymbol{\pi}}
\newcommand{\sangle}{\sphericalangle}
\title{Proofs of Technical Results Justifying an Algorithm of Reactive 3D Navigation of a Mobile Robot through an Unknown Tunnel\footnotemark[1]}
\date{}
\author{Alexey S. Matveev\footnotemark[2] and Andrey V.Savkin\footnotemark[3]}
\begin{document}
\maketitle
\renewcommand{\thefootnote}{\fnsymbol{footnote}}
\footnotetext[1]{This work was supported the Russian Science Foundation 14-21-00041p and by the Australian Research Council.}
\footnotetext[2]{Department of Mathematics and Mechanics, Saint Petersburg State University, St. Petersburg, Russia}
\footnotetext[3]{School of Electrical Engineering and Telecommunications, University of New South Wales, Sydney,  Australia}

\section{Introduction}
 We consider the problem of autonomously driving a mobile robot through an unknown and unstructured 3D tunnel-like environment. This task is elemental for many applications of mobile robotics such as
inspection and servicing of storm runoff networks, channelized aquifers, mines, bypass tunnels for dams, pipelines in sewer networks, power plants, factories, petrochemical, water supply and fluid transportation industries \cite{RohCho02,OnKeBo03,ThBrBe08,WuLiZhShQi08,TurGar10,GaTaTh13,NaPr14},
underwater environmental studies and archeology, exploring flooded cenotes
and ancient cisterns \cite{WHOBGC10full,FaKaJoWe10,MaRiRiCa16full},
indoor and city exploration, rescue, and surveillance by micro air vehicles (MAV) \cite{BaHeRo09,ShBoBeSi09}, to name just a few.
\par
We focus on the case where the robot does not touch the perimeter of the tunnel, should respect a certain safety margin to it, and subject to this limitation, is free to move in all three dimensions. This is typical in aerial and underwater robotics, which are growing fields with ever-extending application possibilities, where they often hold promise
of more effective alternative technology. To acquire this benefit, the robots should be equipped with navigation
algorithms that enable them to operate autonomously through long time and distance scales.
Meanwhile, the studied case does not concern most specialized in-pipe inspection robots manufactured up to now, whose locomotion relies on a firm contact with the tunnel boundary surface \cite{TurGar10,NaPr14}. For them, the proposed navigation solutions substantially benefit from the fact that the workspace is typically structured (is composed of standard element like straight pipe, T-junction, etc.) and may be partly or fully known.
\par
This paper is targeted at another situation of an unstructured and unknown tunnel with a generic local geometry and a robot whose sensory data is confined to a close small patch of tunnel's boundary, whereas the ``opposite patch'' may be out off the sensing distance.
Such confined and overhead unstructured
environments coupled with the lack of an exterior navigation assistance still constitute a real challenge for underwater autonomous navigation so that most respective missions (e.g., in underwater caves or shipwrecks) are performed by human divers or remotely teleoperated robotised vehicles up to now \cite{WHOBGC10full,MaRiRiCa16full}.
Meanwhile, typical approaches to navigation of MAV's\cite{Kendoul12} have
various limitations due to e.g., reliance on known, including pre-installed, landmarks, salient features or environmental patterns \cite{ToVaHoFe06,MoHiKi07,AhLeAnHo09full,LeFrPo11}, maps \cite{TrBuFo05}, a priori image-databases \cite{CoMeGuMa09}, external positioning systems \cite{ScJuRu09full,AhLeAnHo09full}, or a partially structured environment to enable incremental motion planning  \cite{BaPrHeRo11,GrGrBu12full}. Some of these methods, including simultaneous 3D localization and mapping (SLAM), are computationally expensive \cite{AnFiDoMe06}, hardly match capacities of on-board processors and so use extensive off-board computations \cite{BaHeRo09}, which require reliable communication and may cause undesirable feedback delays.
\par
This paper is aimed at showing that
even in the face of nonholonomy, under-actuation, finite control range, poor knowledge of the scene, and limited sensory  capacity, long-distance fully autonomous advancement through a generic tunnel is feasible
at a little computational cost: control signal is generated via a direct reflex-like conversion
of the current observation.
No complex or dubious operations, like building a map or depositing marks in the scene, are employed for navigation.
\par
An extended introduction and discussion of the control law that achieves the stated objective are given in a manuscript submitted by the authors to the IFAC journal {\it Automatica}. The subsequent text basically contains the proofs of the technical facts underlying justification of the convergence and performance of the algorithm proposed in that manuscript, which were not included in it because of the paper length limit.
This text mainly focuses on geometric issues related to the concept of ``tunnel'', its main content may be attributed to the area of differential geometry. The assumptions of the study are borrowed from that manuscript and are partly motivated by navigation and control concerns.
\par
The body of the paper is organized as follows. Sections~\ref{sec1} and \ref{sec.necass} introduce the tunnel environment
and the assumptions of theoretical analysis, respectively.
Section~\ref{sec.recog} presents a particular way of access to local features of tunnel's surface that can be employed for global navigation.
Section~\ref{sec.proper} discusses some properties of tunnels, the main result of Section~\ref{sec.recog} is proved in Section~\ref{app.proofprop}.
\par
Throughout the paper, the following notations are used:
\begin{itemize}
\item $\spr{A}{B}$ and $A\times B$, inner and cross product, respectively;
\item
$\|\cdot\|$, Euclidean norm of a vector and spectral norm of a matrix and linear operator;
\item $\boldsymbol{\pi}(\bldr)$, set of all {\it projections} of point $\bldr \in \br^3$ onto a regular surface $S \subset \br^3$, i.e., points $\bldb \in S$ such that
\begin{equation}
\label{def.dist}
\|\bldr - \bldb\| = \dist{\bldr}{S}:=\min_{\bldb^\prime \in S} \|\bldb^\prime - \bldr\|;
\end{equation}
\item  $\mathfrak{T}_{\bldb}(S)$, plane tangent to $S$ at $\bldb \in S$;
\item  $N(\bldb)$, unit normal to tunnel $S$ directed inside the tunnel;
\item $D_VW$, derivative of the field $W$ in direction of $V$;
\item  $\mathscr{S}_{\bldb}(V) = - D_{V}N$, shape operator (Weingarten map);
\item  $\ff{\bldb}{V}{W}:= \spr{\mathscr{S}_{\bldb}(V)}{W}, V,W \in \ts_{\bldb}(S)$, second fundamental form (shape tensor);
\item  $\varkappa_-(\bldb) \leq \varkappa_+(\bldb)$, principal curvatures, i.e., the eigenvalues of the quadratic form $\ff{\bldb}{\cdot}{\cdot}$;
\item $E_\pm(\bldb)$, respective unit eigenvectors (principal eigenvectors) continuously depending on $\bldb\in S$;
\item  $p_\pm(\bldb)$, lines spanned by $E_\pm(\bldb)$ (principal directions);
\item  $\vec{\tau}(\bldb)\in \ts_{\bldb}(S)$, smooth unit vector-field tangent to the meridian $\mathscr{M}[b(\bldb)]$ that passes through $\bldb$;
\item  $\vec{\tau}_\ast(\bldb_\ast)\in \ts_{\bldb_\ast}[S(d_\ast)]$, similar vector-field on $S(d_\ast)$;
\item $S(d_\ast) := \{\bldp = \bldb + d_\ast N(\bldb), \bldb \in S \setminus \partial S\}$, locus of the desired locations of the robot;
\item  $\mathfrak{T}_{\bldb\ast}[S(d_\ast)]$, plane tangent to $S(d_\ast)$ at $\bldb \in S(d_\ast)$;
\item  $N_\ast(\bldb)$, unit normal to $S(d_\ast)$ directed outward $S$;
\item $\mathbf{Pr}_L$, orthogonal projection onto a subspace $L \subset \br^3$;
\item $\nabla_VW$, Levy-Civita covariant derivative, i.e., $\nabla_VW = \mathbf{Pr}_{\ts} D_V W$, where $\ts$ is the tangent plane;
\item $\mathbf{Id}_{\bldb}$, identity operator acting in $\ts_{\bldb}$;
\item $[A,B,C]:= \spr{A}{B \times C}$, triple product;
\item $\sphericalangle(A,B)$, angle from $A$ to $B$, where $A, B \in \ts_{\bldb}(S)$ and positive angles are countered counterclockwise when looking from the side of $N(\bldb)$;
\item $\mathscr{R}_{\bldb}(\theta)$, rotation of the plane $\ts_{\bldb}(S)$ through angle $\theta$.
\end{itemize}

\section{Tunnel environment}
\label{sec1}
 We consider a mobile robot that operates in a tunnel-like environment bounded by a regular surface $S \subset \br^3$. It is required to advance the robot through the tunnel and, in the case of a long operation, to eventually approach a certain ``comfortable'' value $d_\ast$ of the distance $d$ to $S$.
A troublemaking trait of the robot is that its sensing capacity is limited.
In its local frame, it identifies the direction $\vec{d}$ to a minimum-distance point $\bpi(\bldr)$ of $S$ and the distance to $S$ along any ray emitted from $\bldr$ at an angle $\leq \alpha_s$ with respect to $\vec{d}$, where $\alpha_s>0$ is a constant. Thus the robot ``sees'' only a small patch of $S$ around $\bpi(\bldr)$ 
 but has no access to any global direction, including the desired direction of motion that corresponds to ``advancing along the tunnel''. Thus the robot has no other option but to recognize this direction from features of the sensed small patches of $S$.
\par
 As a result, the control objective is realistic not for any surface $S$: ``advancement along the tunnel'' should be of sense and its direction should be recognizable from observing small patches of $S$. A prototypical sample where these are satisfied is a circular cylinder. Then the ``direction of advancement'' is that of its generatrix, it can be determined as that in which any patch of the boundary has zero curvature. Inspired by this sample, we will deal with surfaces for which contortion of any patch is minimal in an acceptable direction of motion.
 \par
 Now we pass to a rigorous definition.
\par
Let $\boldsymbol{B}$ be either a) the real line $\br$, or b) the unit circle $S^1 \subset \br^2$, or c) an interval $[b_-,b_+] \subset \br$.
\begin{define}
\label{def.tun}
A {\em tunnel with the basis} $\boldsymbol{B}$ is a $C^3$-smooth regular surface $S$ (with a one-dimensional boundary $\partial S = B^{-1}[b_-] \cup B^{-1}[b_+]$ in the case c)) equipped with a {\em projection} $B: S \to \boldsymbol{B}$ such that
\begin{enumerate}[{\bf i)}]
\item The projection is a proper surjective submersion: it is $C^3$-smooth, maps $S$ onto $\boldsymbol{B}$, its differential has rank $1$ everywhere, and the inverse image of any compact subset of $\boldsymbol{B}$ is compact;
\item Any {\em meridian} $\mathscr{M}(b):= B^{-1}[b]$ is a simple closed curve and the set $S$ is closed.
\end{enumerate}
\end{define}
\par
If $\boldsymbol{B} = [b_-,b_+]$, the tunnel is said to be {\it open}, and {\it closed} otherwise. The (multivalued in general) {\it basic coordinate} $b(\bldr):=\{B(\bldb):\bldb \in \bpi(\bldr) \}$ of $\bldr$ assesses progression of robot $\bldr$ through the tunnel: Progression is associated with evolution of $b$ in a certain direction.
The case of $\boldsymbol{B}=\br$ is a model for a ``very long tunnel''.
\par
A prototypical sample of a tunnel is a right cylindrical surface, where $b$ is the coordinate along the axis of symmetry and the meridians are the perpendicular sections.  A closed compact tunnel is exemplified by a torus that is obtained by revolving a circle about an axis coplanar and disjoint with the circle. Then the meridians are instant snapshots of the moving circle and the coordinate $b = (\cos \varphi, \sin \varphi)\in S^1$ assesses the rotation angle $\varphi$. Another example is a surface of revolution, where the generatrix is a graph of a positive function defined on the axis of revolution $AoR$, the meridians are the paths of the points of the revolved curve, and $b$ is the coordinate along $AoR$.
Any diffeomorphism $J$ of an open vicinity of a tunnel $S$ onto an open subset of $\br^3$ transforms $S$ to a new tunnel with the projection $B \circ J^{-1}$ and meridians $J[\mathscr{M}(b)]$. Meanwhile, there exist tunnels of other origins; for example, they may be obtained by properly moving a planar Jordan curve over a space path.
\par
If the tunnel is open, the robot should reach its {\it end}, i.e., $\bpi(\bldr)$ must arrive at the $\delta_s$-vicinity of $\partial S$, where $\delta_s>0$ is given. We assume that the robot can recognize this event via presence of the edge $\partial S$ in the sensory data.
For a closed tunnel, it is needed to ensure that the basic coordinate eventually evolves in an altered direction with a speed separated from zero.
Such a behavior is yet realistic only after tolerating a transient.
\par
With these in mind, we denote by $\bldr(t)$ the position of the robot at time $t$ and arrive at the final statement of the control objective.
\begin{define}
\label{def.evadv}
The robot is said to {\em solve the tunnel} if
$
d(t) := \dist{\bldr(t)}{S} = \|\vec{d}(t)\| >0 \; \forall t
$ and there exists time $t_0$ such that
\begin{itemize}
\item For an open tunnel, $\bpi[\bldr(t_0)]$ is in the $\delta_s$-vicinity of $\partial S$;
\item In the case of a closed tunnel, the following statements are true whenever $t \geq t_0$:
\begin{enumerate}[{\bf i)}]
\item The robot's basic coordinate $b(t)$ is unique, smoothly depends on time, and $\pm \dot{b}(t) \geq v_b$, where the sign in $\pm$ and $v_b>0$ do not alter with time;
\item The distance $d(t)$ to $S$ smoothly depends on time and monotonically converges to $d_\ast$ as $t \to \infty$.
\end{enumerate}
\end{itemize}
\end{define}
\par
This covers the cases of operation inside the tunnel and outside it.
For the sake of definiteness, we focus on the first case.
If $\boldsymbol{B}= S^1$, the derivative $\dot{b}$ is meant as that of the angular coordinate of $b \in S^1$.
\par
If the tunnel is open, Definition~\ref{def.evadv} sets an objective to run the entire length of the tunnel from an initial location to some end of the tunnel and, thus, to arrive at a position of leaving the tunnel. If the tunnel is closed, an endless progression through the tunnel in a fixed direction is targeted, with the distance $d$ approaching $d_\ast$.
\par
Now we introduce a class of tunnels for which the direction of motion is recognizable from any its small patch.
\begin{define}
\label{def.locdir}
A tunnel is said to be {\em regular} if for any $\bldb \in S$, the normal curvature of $S$ in the meridian direction $\vec{\tau}(\bldb)$ is not minimal, moreover, it is separated from $\varkappa_-(\bldb)$ by a gap $\Delta_{\tau}>0$ that does not depend on $\bldb$:
\begin{equation}
\label{basic.ineq}
\ff{\bldb}{\vec{\tau}(\bldb)}{\vec{\tau}(\bldb)} \geq \varkappa_-(\bldb) + \Delta_\tau .
\end{equation}
\end{define}
The second claim follows from the first one if $S$ is compact. Due to \eqref{basic.ineq}, $\varkappa_-(\bldb) < \varkappa_+(\bldb)$ and so the principal directions $p_\pm(\bldb)$ and vector-fields $E_\pm$ are well-defined.
\par
The right cylindrical surface is a regular tunnel since $\varkappa_-=0$ and the meridian curvatures are positive.
For the torus, the curvature of any meridian is $r^{-1}$ and $\varkappa_-=\frac{\cos u}{R+r\cos u}$ \cite[p.~157]{Carmo76}, where $R>r$ and the torus is parametrized by $[(R+r\cos u)\cos v, (R+r\cos u)\sin v, r \sin u], u,v \in [0, 2\pi)$. Hence the torus is a regular tunnel since $ r^{-1}>\varkappa_-$.
\par
From now on, we deal with a regular tunnel.
\par
In conclusion of the section, we discover a geometric sense of i) in Definition~\ref{def.evadv} by using
the map
\begin{equation}
\label{map.j}
\bldb \in S \setminus \partial S \mapsto J(\bldb):=\bldb + d_\ast N(\bldb) \in S(d_\ast),
\end{equation}
which diffeomorphically maps $S \setminus \partial S$ onto $S(d_\ast)$ (under the assumptions of this paper, see Corollary~\ref{cor.smooth}). This map gives rise to the {\it meridians $\mathscr{M}_\ast(b) = J[\mathscr{M}(b)]$ on} $S(d_\ast)$. Their tangents can be arranged into the smooth unit vector-field $\vec{\tau}_\ast$ on $S(d_\ast)$:
\begin{equation}
\label{tau.ast}
\vec{\tau}_\ast(\bldr) := \frac{J^\prime(\bldb) \vec{\tau}(\bldb)}{\|J^\prime(\bldb) \vec{\tau}(\bldb)\|}, \quad \text{where} \; \bldb := J^{-1}(\bldr).
\end{equation}
\par
In Definition~\ref{def.evadv}, i) means that while moving nearly over $S(d_\ast)$, the robot transverses the meridians in an unaltered direction at angles that are separated from zero and eventually crosses any meridian that lies in the direction of motion with respect to the initial one.

\section{Assumptions of theoretical analysis}
\label{sec.necass}
\setcounter{equation}{0}
\begin{assumption}
\label{ass.unique}
For any $\bldr \in S(d_\ast)$, the set $\bpi(\bldr)$ of projections of $\bldr$ onto $S$ contains a unique point $\bldb(\bldr)$.
\end{assumption}
For convenience's sake, we define the operational zone $Z_{\text{op}}$ of the robot in terms of the extreme values $d_- < d_+$ that are assumed by the distance $d$ to the surface $S$ in this zone:
\begin{equation}
\label{dop}
Z_{\text{op}} := \{\bldr = \bldb + d N(\bldb): \bldb \in S \setminus \partial S, d \in (d_-,d_+)\},
\qquad
\text{where} \quad 0 < d_- < d_\ast < d_+.
\end{equation}
\begin{assumption}
\label{as.cum}
Assumption~{\rm \ref{ass.unique}} is valid with any $d_\ast$ from $[0,d_+]$ and there exists $\Delta_\varkappa \in (0,1]$ such that
\begin{equation}
\label{ineq.cur2}
 1- d_+ \varkappa_+(\bldb) \geq \Delta_\varkappa \qquad \forall \bldb \in S.
\end{equation}
\end{assumption}
Whenever $\varkappa_+(\bldb) \leq 0$, \eqref{ineq.cur2} is certainly true.
\par
The last assumption is purely technical.
\begin{assumption}
\label{ass.lipsch}
The maps $N, \varkappa_\pm, E_\pm,\nabla B$ are Lipschitz:
there exist constants $L_N,L_\varkappa, L_E, L_B>0$ such that
\begin{gather}
\label{lip.const}
\| N(\bldb_1)-  N(\bldb_2)\| \leq L_N \|\bldb_1- \bldb_2\|,
\\
\nonumber
| \varkappa_\pm(\bldb_1)-  \varkappa_\pm (\bldb_2)| \leq L_\varkappa \|\bldb_1- \bldb_2\|,
\\
\nonumber
\| E_\pm(\bldb_1)-  E_\pm (\bldb_2)\| \leq L_E \|\bldb_1- \bldb_2\|,
\\
\nonumber
\| \nabla B(\bldb_1)-  \nabla B(\bldb_2)\| \leq L_B \|\bldb_1- \bldb_2\| \; \forall \bldb_1, \bldb_2 \in S.
\end{gather}
Furthermore, there exist $\Delta_B^\pm>0$ such that $\Delta_B^- \leq \|\nabla B(\bldb)\| \leq \Delta_B^+ \; \forall \bldb \in S$.
\end{assumption}
\par
For a compact tunnel, this assumption necessarily holds.
\section{Direction estimator}
\label{sec.recog}
\setcounter{equation}{0}
Its role is to generate the direction of motion along the tunnel: This direction should be transversal to the meridian passing through $\bldb$ to meet i) in Definition~\ref{def.evadv}. The premise is the sensory data, i.e., the direction $\vec{d}$ from robot's location $\bldr$ to its projection $\bldb:=\bpi(\bldr)$ onto $S$ and the distance $d(\alpha,\varphi)$ to $S$ along any ray $\mathfrak{R}$ emitted from $\bldr$ at an angle $\alpha \in [0, \alpha_s]$ to $\vec{d}$. Here $\varphi$ is the angle giving the direction of $\mathfrak{R}$ in projection onto the plane $\vec{d}^{\,\bot}$ normal to $\vec{d}$.
\par
As a basic choice, the principal direction $p_-(\bldb)$ might be used since it transverses the respective meridian due to \eqref{basic.ineq}.
However, computation of $p_-(\bldb)$ involves second-order differentiation of the sensory data, which is a highly unstable procedure. So we employ differentiation-free and simpler methods. We do not confine ourselves to $p_-(\bldb)$ and a particular method but adopt a whole class of them that is delineated in the following definition. It takes into account that the direction of motion in needed only in the {\it active operational zone} $Z_{\text{op}}^a :=\{\bldr \in Z_{\text{op}}: \bpi(\bldr)\,\text{is not in the $\delta_s$-vicinity of} \, \partial S\}$ since the mission is terminated as soon as $\bpi(\bldr)$ enters this vicinity. For closed tunnels, $\partial S = \emptyset$ and so $Z_{\text{op}}^a =Z_{\text{op}} $.
\begin{define}
\label{dir.estdef}
{\em Direction estimator (DE)} maps location dependent $O=O(\bldr)$ sensory data into a straight line $p(\bldr) \subset \ts_{\bpi(\bldr)}(S)$ such that the associated map $\bldr \in Z_{\text{\rm op}}^a \mapsto p(\bldr) \in \mathbf{Gr}(1,3)$ is continuous. Its {\em exactness} is an upper estimate $\beta>0$ of the angle between the lines $p(\bldr)$ and $p_-[\bpi(\bldr)]$ that holds for any $\bldr \in Z_{\text{\rm op}}^a$.
\end{define}
Here $\mathbf{Gr}(1,3)$ is the smooth Grassmanian manifold of all one-dimensional linear subspaces of $\br^3$ \cite[pp.~42-44]{Shaf94}.
We assume that the output $p(O)$ of DE is given in the local frame of the robot.
For a DE to be useful, its exactness $\beta$ should be less than the minimal angular discrepancy between $p_-[\bpi(\bldr)]$ and the meridian direction at point $\bpi(\bldr)$ since then $p(\bldr)$ is transversal to the meridian, as is required.
\par
Now we discuss a particular design of DE under which arbitrarily high exactness can be achieved.
\par
{\bf Most-distant-point-based estimator} (MDPBE) with parameter $\alpha_{e} \in (0,\alpha_s]$ finds the tangential directions (given by $\varphi$) of the ray $\mathfrak{R}$ that are
associated with the local maxima of the distance $d(\alpha_e,\varphi)$, shifts every of the found $\varphi$'s into $(-\pi/2,\pi/2]$ by adding, if necessary, an integer multiple of $\pi$, computes the arithmetic mean $\varphi_{\to}$ of the resultant angles, and returns the straight line $p$ that goes in the normal plane $\vec{d}^{\,\bot}$ in the direction of $\varphi_{\to}$. MDPBE is well-posed if the number of local maxima is finite.
\begin{proposition}
\label{prop.de}
Let Assumptions~{\rm \ref{as.cum}} and {\rm \ref{ass.lipsch}} hold.
For any $\beta>0$, there is $\ov{\alpha} \in (0,\alpha_s]$ such that
MDPBE is well-posed (with two local maxima) in the active operational zone and is a direction estimator with exactness $\beta$ whenever $\alpha_e \in (0,\ov{\alpha})$.
\end{proposition}
\par
The proof of this proposition is given in Section~\ref{app.proofprop}.
\section{General properties of tunnels and technical facts}
\label{sec.proper}
Let Assumptions~\ref{as.cum}---\ref{lip.const} hold and let $L_N$ be taken from \eqref{lip.const}.
\begin{lemma}
\label{lem.cuvb}
The norm of the operator $\mathscr{S}_{\bldb}: \ts_{\bldb}(S) \to \ts_{\bldb}(S)$ and the principal curvatures obey the estimates
\begin{equation}
  \label{est.norm}
  \|\mathscr{S}_{\bldb}\| \leq L_N \quad \text{\rm and} \quad |\varkappa_\pm(\bldb)| \leq L_N \quad \forall \bldb \in S .
\end{equation}
There exists $\ve>0$ such that $(d,\bldb) \mapsto h(d,\bldb):=\bldb+d N(\bldb)$ is a $C^2$-diffeomorphism of $T:=(-\ve,d_+) \times (S\setminus \partial S)$ onto an open neighborhood of \eqref{dop}, and $h^{-1}(\bldp) = \{\dist{
\bldp}{S},\bpi(\bldp)\} , \bpi(\bldp) \not \in \partial S \; \forall \bldp \in Z_{\text{\rm op}}$.
\end{lemma}
\pf
Inequalities \eqref{est.norm} hold since $\|\mathscr{S}_{\bldb} V\|=\|D_V N\| \leq L_N \|V\|$ for all $V \in \ts_{\bldb}(S)$ by \eqref{lip.const} and $\varkappa_\pm$ are the eigenvalues of $\mathscr{S}_{\bldb}$.
By \eqref{ineq.cur2} and \eqref{est.norm}, there exists $\ve\in (0, d_+/2)$ such that $1-d \varkappa_+(\bldb) \geq \Delta_\varkappa \; \forall \bldb \in S, d \in [0,d_+]$ and $1-d \varkappa_-(\bldb) \geq \Delta_\varkappa \; \forall \bldb \in S, d \in [-\ve,0]$. Hence for $d \in [-\ve, d_+]$ and $\bldb \in S$, the symmetric operator $\mathbf{Id}_{\bldb} - d\cdot \mathscr{S}_{\bldb}: \ts_{\bldb}(S) \to \ts_{\bldb}(S)$ is positively definite, its least eigenvalue is no less than $\Delta_\varkappa$ and so
\begin{equation}
\label{norm.inv}
\left\| (\mathbf{Id}_{\bldb} - d\cdot \mathscr{S}_{\bldb})^{-1} \right\| \leq \Delta_\varkappa^{-1}.
\end{equation}
The differential of $h$ maps $(\eta,V) \in \br \times \ts_{\bldb}(S)$ into $ W =  \eta N + (\mathbf{Id}_{\bldb} - d\mathscr{S}_{\bldb}) V \in \br^3$ and is invertible:
\begin{equation}
\label{form.inv}
\eta = \spr{W}{N}, \quad V = (\mathbf{Id} - d\mathscr{S}_{\bldb})^{-1} [W - \eta N].
\end{equation}
Hence $h$ is a local diffeomorphism. So the image $h(T)$ is open and the inverse image $h^{-1}(\bldp)$ is a discrete set for any point $\bldp \in h(T)$. For any $(d,\bldb) \in h^{-1}(\bldp)$, we have $\bldp = h(d,\bldb)= \bldb + d N(\bldb), |d| \leq d_+ \Rightarrow \|\bldb\| \leq \|\bldp\| + d_+$ and so the inverse image $h^{-1}(\bldp)$ is bounded.
Suppose that it is infinite. Then there exists an infinite sequence $\{(d_j,\bldb_j)\} \subset h^{-1}(\bldp)$ such that $d_j \to d, \bldb_j \to \bldb$ as $j\to \infty$ and $\zeta_j:=(d_j,\bldb_j) \neq \zeta:=(d,\bldb)\; \forall j$. Here $d \in [-\ve,d_+]$ and $\bldb \in S$ and so the differential $h^\prime(d,\bldb)$ is invertible by the foregoing. Meanwhile for any limit point $V$ of the bounded sequence $(\zeta_j-\zeta)/\|\zeta_j-\zeta\|$, we have $\|V\| =1$ and $h^\prime(d,\bldb) V=0 \Rightarrow V=0 \Rightarrow \|V\|=0$. The contradiction obtained proves that the set $h^{-1}(\bldp)$ is finite.
\par
This and local diffeomorphic property of $h(\cdot)$ imply that
$h(\cdot)$ is a covering map \cite[p.~67]{BER88}.
Definition~\ref{def.tun} and Assumption~\ref{ass.lipsch} imply that $S, S \setminus \partial S$ and so $T$ are arcwise connected. Hence all fibres $F_{\bldp}:= h^{-1}(\bldp), \bldp \in h(T)$ have the same size $M$ (the multiplicity of the covering)
\cite[Th.~2.4.4]{BER88}. We are going to show that $M=1$, possibly after properly decreasing $\ve$.
\par
To this end, we first show this for the restriction $h|_{(-\ve,\ve) \times (S\setminus \partial S)}$, which evidently also is a local diffeomorphism and a covering map. Suppose to the contrary that its multiplicity is greater than $1$ for $\ve = \ve_j \; \forall j$, where $\ve_j\in (0, \ve)$ and $\ve_j \to 0$ as $j \to \infty$. We also pick $\bldb \in S\setminus \partial S$. For any $j$, the fibre $F_{\bldb}$ contains $(0,\bldb)$ and some other point $(\eta_j,\bldb_j) \in (-\ve_j,\ve_j)\times (S \setminus \partial S)$, where $\bldb = h(0,\bldb) = h(\eta_j,\bldb_j) = \bldb_j + \eta_j N[\bldb_j]$ and so $\eta_j \neq 0$. Hence
$$
\|\bldb_j - \bldb\| = \eta_j \leq \ve_j, \quad \spr{\eta_j^{-1}(\bldb-\bldb_j)}{N(\bldb)} =  \spr{ N[\bldb_j]}{N[\bldb]}.
$$
As $j \to \infty$, the unit vector $\eta_j^{-1}(\bldb-\bldb_j)$ converges to the tangent plane $\mathfrak{T}_{\bldb}(S)$ and so the left hand side of the last equation goes to zero, whereas its right hand side goes to $1$. This contradiction proves that for all small enough $\ve>0$, the multiplicity of the map $h|_{(-\ve,\ve) \times (S\setminus \partial S)}$ is $1$ and so this map is a diffeomorphism onto its image. We assume that $\ve$ is decreased (if necessary) to this range.
\par
Now suppose that that the multiplicity $M$ of $h(\cdot)$ on $T$ is greater than $1$, pick $\bldb_1 \in S\setminus \partial S \subset h(T)$, and enumerate $(d_1=0,\bldb_1), (d_2,\bldb_2), \ldots, (d_M,\bldb_M)$ the elements of the fibre $F_{\bldp}$, where $\bldp := h(0,\bldb_1) = \bldb_1$. Here
\begin{equation}
\label{rome}
\bldb_1 = h(d_i, \bldb_i) = \bldb_i+d_i N(\bldb_i)
\end{equation}
and $d_i \in [\ve ,d_+)\; \forall i \geq 2$ since $d_i \in (-\ve, \ve)$ is excluded by the just established property of $h|_{(-\ve,\ve) \times (S\setminus \partial S)}$. For $i \geq 2$, the straight-line segment connecting $\bldb_i$ and $\bldb_1$ is a subset of $h([0,d_i]\times \{\bldb_i\}) \subset h(T)$ and so
 the point $(d_i,\bldb_i) \in F_{\bldp}$ gives rise to a lift of the path $\mathfrak{p}(\theta) = (1-\theta)\bldb_1 + \theta \bldb_2 \in h(T), \theta \in [0,1]$ that starts from this point, i.e., to a continuous map $L_i(\cdot) = [d_i(\cdot), \bldb_i(\cdot)]: [0,1] \to T$ such that
\begin{equation}
\label{lift}
h[L_i(\theta)] = \mathfrak{p}(\theta)\qquad \forall \theta
\end{equation}
and $[d_i(0),\bldb_i(0)] = [d_i,\bldb_i]\; \forall i$. Since the path $\mathfrak{p}(\cdot)$ is smooth and $h(\cdot)$ is locally diffeomorphic, any lift $L_i(\cdot)$ is also smooth. For $\theta \in [0,1]$, the points $L_1(\theta), \ldots, L_M(\theta)$ are pair-wise distinct and exhaust the fibre $h^{-1}[\mathfrak{p}(\theta)] = F_{\mathfrak{p}(\theta)}$.
 \par
By Assumption~\ref{as.cum} (in the part that refers to Assumption~\ref{ass.unique}), the projection $\bldb(\theta)$ of $\mathfrak{p}(\theta)$ onto $S$ is unique for all $\theta \in [0,1]$. So not only $\dist{\mathfrak{p}(\theta)}{S}$ but also $\bldb(\theta)$ continuously depend on $\theta$. Furthermore,
\begin{equation}
\label{norm.dist}
\bldb(0) = \bldb_1, \quad \bldb(1) = \bldb_2, \qquad \text{and} \quad
\mathfrak{p}(\theta) = \bldb(\theta) \pm \dist{\mathfrak{p}(\theta)}{S} N[\bldb(\theta)] \quad \text{whenever} \; \bldb(\theta) \not\in \partial S.
\end{equation}
Since $0\leq \dist{\mathfrak{p}(\theta)}{S}<\ve$, $-\ve< d_1(\theta) < \ve$, and $\bldb(\theta) \approx \bldb_1 \in S \setminus \partial S \Rightarrow \bldb(\theta) \in S \setminus \partial S$ for $\theta \approx 0$, the diffeomorphic property of $h|_{(-\ve,\ve)\times (S\setminus \partial S)}$ and \eqref{lift} (with $i:=1$) yield that
\begin{equation}
\label{fhglg}
d_1(\theta) = \dist{\mathfrak{p}(\theta)}{S}, \quad \bldb_1(\theta) = \bldb(\theta) \qquad \forall \theta \in [0,\tau]
\end{equation}
with sufficiently small $\tau \in (0,1]$. We are going to show that in fact this is true with $\tau=1$.
\par
Suppose to the contrary that the maximal such $\tau$ is less than $1$. Then $[\dist{\mathfrak{p}(\theta)}{S}, \bldb(\theta)] \neq [d_1(\theta), \bldb_1(\theta)]$ for $\theta = \theta_j \; \forall j$, where $\theta_j \in (\tau,1)$ and $\theta_j \to \tau$ as $j \to \infty$. So \eqref{lift}, \eqref{norm.dist} and the diffeomorphic property of $h|_{(-\ve,\ve)\times (S\setminus \partial S)}$ exclude the case where $\dist{\mathfrak{p}(\tau)}{S}=0$. Thus
\begin{equation}
\label{positive}
d_1(\tau)=\dist{\mathfrak{p}(\tau)}{S} >0
\end{equation} and so
$\dist{\mathfrak{p}(\theta_j)}{S}, d_i(\theta_j) \in (0, d_+)$ and $\bldb(\theta_j) \not\in \partial S$ for $j \approx \infty$. From now on, we focus on only such $j$'s.
\par
Suppose that there are arbitrarily large $j$'s for which \eqref{norm.dist} with $\theta:=\theta_j$ holds with the sign $+$. Then passing to a proper subsequence ensures that this is true for all $j$. In this case,  $[\dist{\mathfrak{p}(\theta_j)}{S}, \bldb(\theta_j)] \in T, h [\dist{\mathfrak{p}(\theta_j)}{S}, \bldb(\theta_j)] = \mathfrak{p}(\theta_j)$, i.e., $[\dist{\mathfrak{p}(\theta_j)}{S}, \bldb(\theta_j)] \in F_{\mathfrak{p}(\theta_j)}$. Hence
$[\dist{\mathfrak{p}(\theta_j)}{S}, \bldb(\theta_j)] = [d_i(\theta_j),\bldb_i(\theta_j)]$ for some $i=i(j) =1,\ldots,M$. Here $i(j) \neq 1$ since
$[\dist{\mathfrak{p}(\theta_j)}{S}, \bldb(\theta_j)] = [d_1(\theta_j),\bldb_1(\theta_j)]$.
Then by passing to a subsequence once more, it can be ensured that $i(j)$ does not alter with $j$, i.e., $i(j) \equiv i \geq 2 \; \forall j$. Then passing to the limit as $j\to \infty$ in the equation $\dist{\mathfrak{p}(\theta_j)}{S} = d_i(\theta_j)$ yields $\dist{\mathfrak{p}(\tau)}{S} = d_i(\tau)$. It follows that the distance $\dist{\mathfrak{p}(\tau)}{S}$ from the point $\mathfrak{p}(\tau)$ to the surface $S$ is attained at both $\bldb_1(\tau)$ and $\bldb_i(\tau)$, where $[d_1(\tau), \bldb_1(\tau)] \neq [d_i(\tau), \bldb_i(\tau)]$ and $d_1(\tau) = d_i(\tau) = \dist{\mathfrak{p}(\tau)}{S}$ imply that  $\bldb_1(\tau) \neq \bldb_i(\tau)$, in violation of Assumption~\ref{as.cum}. The contradiction obtained proves that \eqref{norm.dist} with $\theta:=\theta_j$ holds with the sign $-$ for all large enough $j$'s. Then by \eqref{lift} and \eqref{norm.dist}, we have
\begin{gather*}
\bldb(\theta_j) - \dist{\mathfrak{p}(\theta_j)}{S} N[\bldb(\theta_j)] = \mathfrak{p}(\theta_j) =\bldb_1(\theta_j) + d_1(\theta_j) N[\bldb_1(\theta_j)] \\
\Downarrow j \to \infty
\\
\bldb(\tau) - \dist{\mathfrak{p}(\tau)}{S} N[\bldb(\tau)] = \mathfrak{p}(\tau) =\bldb_1(\tau) + d_1(\tau) N[\bldb_1(\tau)]
\overset{\text{\eqref{fhglg}}}{\Rightarrow} \dist{\mathfrak{p}(\tau)}{S} =0,
\end{gather*}
in violation of \eqref{positive}. The contradiction obtained proves that \eqref{fhglg} holds with $\tau=1$.
\par
It follows that $d_1(1) =0 \Rightarrow \bldb_2= \mathfrak{p}(1) = \bldb_1(1) + d_1(1) N [\bldb_1(1)] = \bldb_1(1)$. Hence
\begin{gather*}
\bldb_1(\theta) + d_1(\theta) N[\bldb_1(\theta)] = h[d_1(\theta), \bldb_1(\theta)] = \mathfrak{p}(\theta) = (1-\theta)\bldb_1 + \theta \bldb_2
\overset{\text{\eqref{rome}}}{=} (1-\theta)[\bldb_2+d_2 N(\bldb_2)] + \theta \bldb_2
\\
= \bldb_2 +d_2(1-\theta) N (\bldb_2) \qquad \forall \theta \in [0,1]
\\
\Downarrow
\\
d_1^\prime(\theta) N[\bldb_1(\theta)] + \left[ \mathbf{Id}_{\bldb_1(\theta)} - d_1(\theta) \mathscr{S}_{\bldb_1(\theta)} \right] \bldb_1^\prime(\theta) 
= \frac{d}{d \theta} \left\{ \bldb_2 +d_2(1-\theta) N (\bldb_2) \right\} = -d_2 N [\bldb_1(1)] .
\end{gather*}
By putting $\theta:= 1$ here and taking into account that both $\bldb_1^\prime(\theta)$ and $\left[ \mathbf{Id}_{\bldb_1(\theta)} - d_1(\theta) \mathscr{S}_{\bldb_1(\theta)} \right] \bldb_1^\prime(\theta)$ lie in the tangential plane, whereas $N[\bldb_1(\theta)]$ is normal to it, we see that
$d_2=-d_1^\prime(1)$. Thus
\begin{gather*}
\left[ \mathbf{Id}_{\bldb_1(\theta)} - d_1(\theta) \mathscr{S}_{\bldb_1(\theta)} \right] \bldb_1^\prime(\theta) = d_1^\prime(1) N [\bldb_1(1)] - d_1^\prime(\theta) N[\bldb_1(\theta)]
\\
=
d_1^\prime(1) \left\{ N [\bldb_1(1)]- N[\bldb_1(\theta)] \right\} - [d_1^\prime(\theta) - d_1^\prime(1)] N[\bldb_1(\theta)].
\end{gather*}
Then after applying the orthogonal projection onto the tangent plane $\mathfrak{T}_{\bldb_1(\theta)}(S)$, we see that
$$
\left\| \left[ \mathbf{Id}_{\bldb_1(\theta)} - d_1(\theta) \mathscr{S}_{\bldb_1(\theta)} \right] \bldb_1^\prime(\theta) \right\|
\leq d_2 \left\| N [\bldb_1(1)]- N[\bldb_1(\theta)] \right\| \overset{\text{\eqref{lip.const}}}{\leq} L_N d_2 \left\| \bldb_1(1)- \bldb_1(\theta) \right\|.
$$
On the other hand, \eqref{norm.inv} implies that the left hand side is no less than $ \| \bldb_1^\prime(\theta)\|/\Delta_\varkappa $.
As a result,
$$
\| \bldb_1^\prime(\theta)\| \leq L_N \Delta_\varkappa d_2 \left\| \bldb_1(1)- \bldb_1(\theta) \right\|
$$
and so $\bldb_1(\theta) = \bldb_1(1)\; \forall \theta$ by Gr\"{o}wall's inequality. Putting $\theta:= 0$ here yields $\bldb_1=\bldb_1(0) = \bldb_1(1) = \bldb_2$, in violation of the foregoing. This contradiction proves that the multiplicity $M=1$ and so $h(\cdot)$ is a diffeomorphism.
\par
To prove the last claim of the lemma, we consider $\bldp \in Z_{\text{op}}$. By \eqref{dop}, $\bldp = \bldb + d N(\bldb)$, where $\bldb \in S \setminus \partial S, d \in (d_-,d_+)$. We put
\begin{equation}
\label{defpnew}
\mathfrak{p}(\theta) := \bldb + d \theta N(\bldb)\qquad \theta \in [0,1]
 \end{equation}
 and partly retrace the foregoing arguments. Specifically, we first establish that $\bldb(\theta) = \bldb$ and $\dist{\mathfrak{p}(\theta)}{S}=d \theta$ for all small enough $\theta$ by using the diffeomorphic property of the restriction $h|_{(-\ve,\ve) \times (S\setminus \partial S)}$. As before, we then introduce the maximal interval $[0,\tau], \tau \leq 1$ on which \eqref{fhglg} holds with
$d_1(\theta) := \theta d$ and $\bldb_1(\theta) := \bldb$:
\begin{equation}
\label{dfkg}
\theta d = \dist{\mathfrak{p}(\theta)}{S}, \quad \bldb = \bldb(\theta) \qquad \forall \theta \in [0,\tau].
\end{equation}
It remains to show that $\tau=1$. Suppose to the contrary that $\tau <1$. Then there there exists an infinite sequence $\{\theta_j\} \subset (\tau,1]$ such that
\begin{equation}
 \label{nefh}
 [\dist{\mathfrak{p}(\theta_j)}{S}, \bldb(\theta_j)] \neq [\theta_j d, \bldb], \quad \dist{\mathfrak{p}(\theta_j)}{S} \in [0,d_+), \quad \bldb(\theta_j) \not\in \partial S \qquad \forall j
 \end{equation}
 and $\theta_j \to \tau$ as $j \to \infty$. By passing to a subsequence, if necessary, it can be ensured that \eqref{norm.dist} with $\theta:=\theta_j$ holds either with the sign $+$ for all $j$ or with the sign $-$ for all $j$. However in the first case, we have $h [\dist{\mathfrak{p}(\theta_j)}{S}, \bldb(\theta_j)]= h[\theta_j d, \bldb]$, in violation of the first inequality from \eqref{nefh} thanks to the established diffeomorphic property of $h(\cdot)$. So the second option holds: $\mathfrak{p}(\theta_j) = \bldb(\theta_j) - \dist{\mathfrak{p}(\theta_j)}{S} N[\bldb(\theta_j)]$. Letting $j \to \infty$ yields that
 $$
 \bldb + d \tau N(\bldb) \overset{\eqref{defpnew}}{=}\mathfrak{p}(\tau) = \bldb(\tau) - \dist{\mathfrak{p}(\tau)}{S} N[\bldb(\tau)] \overset{\text{\eqref{dfkg}}}{=} \bldb - \tau d N[\bldb] \Rightarrow d =0,
 $$
in violation of $d \in (d_-,d_+)$. This contradiction completes the proof.
\epf
\begin{corr}
\label{cor.smooth}
The map \eqref{map.j} is diffeomorphic, the norms of the differentials of $J$ and $J_{-1}:=J^{-1}$ do not exceed $1+d_\ast L_N$ and  $1/\Delta_\varkappa$, respectively, and  $S(d_\ast)$ is a smooth surface.
\end{corr}
\begin{lemma}
\label{lem.angles}
Let $Q : \ts_{\bldb}(S) \to \ts_{\bldb}(S) $ be a linear symmetric positively definite operator with the eigenvalues $0 < q_- \leq q_+ $. Then
$$
\sin \sphericalangle(QA,QB) = \zeta \sin \sphericalangle(A,B), \; \text{\rm where}\; \zeta \in \left[ q_-/q_+, q_+/q_-\right].
$$
\end{lemma}
\pf For all tangent vectors $C,D \in \ts_{\bldb}(S)$, we have
\begin{equation}
\label{sin.angle}
\sin \sphericalangle(C,D)=[N,C,D]/(\|C\|\|D\|) .
\end{equation}
Let $\e_-$ and $\e_+$ be orthonormal eigenvectors of $Q$ associated with $q_-$ and $q_+$, respectively. Since $QC = q_-\spr{C}{\e_-} \e_- + q_+\spr{C}{\e_+}\e_+$, we have
\begin{gather*}
\sin \sphericalangle(QA,QB) =  \frac{[N,QA,QB]}{\|QA\| \|QB\|}
\\
= \frac{\spr{N}{QA \times QB}}{\|QA\| \|QB\|} = \frac{\spr{N}{[q_-\spr{A}{\e_-} \e_- + q_+\spr{A}{\e_+}\e_+] \times [q_-\spr{B}{\e_-} \e_- + q_+\spr{B}{\e_+}\e_+]}}{\|QA\| \|QB\|}
\\
=
\frac{\spr{N}{[q_-\spr{A}{\e_-} \e_-] \times [q_-\spr{B}{\e_-} \e_-]}}{\|QA\| \|QB\|}
+
\frac{\spr{N}{[q_-\spr{A}{\e_-} \e_- ] \times [q_+\spr{B}{\e_+}\e_+]}}{\|QA\| \|QB\|}
\\
+
\frac{\spr{N}{[q_+\spr{A}{\e_+}\e_+] \times [q_-\spr{B}{\e_-} \e_- ]}}{\|QA\| \|QB\|}
+
\frac{\spr{N}{[q_+\spr{A}{\e_+}\e_+] \times [q_+\spr{B}{\e_+}\e_+]}}{\|QA\| \|QB\|}
\\
= q_-q_+\frac{\spr{N}{[\spr{A}{\e_-} \e_- ] \times [\spr{B}{\e_+}\e_+]}}{\|QA\| \|QB\|} + q_-q_+\frac{\spr{N}{[\spr{A}{\e_+}\e_+] \times [\spr{B}{\e_-} \e_- ]}}{\|QA\| \|QB\|}
\\
= q_-q_+\frac{\spr{N}{[\spr{A}{\e_-} \e_- + \spr{A}{\e_+} \e_+] \times [\spr{B}{\e_+}\e_+]}}{\|QA\| \|QB\|} + q_-q_+\frac{\spr{N}{[\spr{A}{\e_-}\e_-+ \spr{A}{\e_+}\e_+] \times [\spr{B}{\e_-} \e_- ]}}{\|QA\| \|QB\|}
\\
= q_-q_+\frac{\spr{N}{[A \times [\spr{B}{\e_-}\e_-+\spr{B}{\e_+}\e_+]}}{\|QA\| \|QB\|} = q_-q_+\frac{\spr{N}{A \times B}}{\|QA\| \|QB\|}
\\
=q_-q_+ \frac{[N,A,B]}{\|QA\| \|QB\|} = q_-q_+ \frac{\sin \sphericalangle(A,B)}{\|QA\| \|QB\|} \|A\|\|B\|.
\end{gather*}
It remains to note that $q_- \|C\|\leq \|QC\| \leq q_+\|C\|$. \epf
\par
For any two vectors $A,B$, we denote by $B_{A^\bot} := B - \spr{B}{A}A/\|A\|^2$ the orthogonal projection of $B$ onto the plane normal to $A$.
We also note that if $A,B \in \ts_{\bldr\ast}[S(d_\ast)] (\bldr \in S(d_\ast))$, then $B_{A^\bot} \in \ts_{\bldr\ast}[S(d_\ast)]$ is the orthogonal projection in the tangent plane $\ts_{\bldr\ast}[S(d_\ast)]$ of vector $B$ onto the line normal to $A$.
\begin{lemma}
\label{lem.rotat}
Let a point $\bldr = \bldr(t)$ smoothly move over $S(d_\ast)$ and let $V(t),W(t) \in \ts_{\bldr(t)\ast}[S(d_\ast)]$ be smooth non-vanishing vector-fields defined on its trajectory. For $\phi:= \sangle (V,W)$, we have
\begin{equation*}
\dot{\phi}\cos \phi = \frac{\spr{(\nabla_{\dot{\bldr}} V)_{V^\bot}}{\mathscr{R}(-\frac{\pi}{2})W} - \spr{(\nabla_{\dot{\bldr}} W)_{W^\bot}}{\mathscr{R}(-\frac{\pi}{2})V}}{\|V\|\|W\|}.
\end{equation*}
\end{lemma}
\pf
We put $V_0=V/\|V\|, W_0=W/\|W\|$.
By \eqref{sin.angle},
\begin{gather*}
 \dot{\phi}\cos \phi = \frac{d}{dt}[N_\ast,V_0,W_0]
= [\dot{N}_\ast,V_0,W_0] + [N_\ast,\dot{V}_0,W_0] + [N_\ast,V_0,\dot{W}_0].
\end{gather*}
Here $\|N_\ast\| \equiv 1 \Rightarrow \spr{\dot{N}_\ast}{N_\ast}=0 \Rightarrow \dot{N}_\ast\in \ts_\ast[S(d_\ast)]$.
 Hence the three vectors $\dot{N}_\ast,V_0,W_0$ lie in a common (tangent) plane. Thus they are linearly dependent and so $[\dot{N}_\ast,V_0,W_0]=0$. Meanwhile,
 $$
 \nabla_{\dot{\bldr}} V = \mathbf{Pr}_{\ts_\ast[S(d_\ast)]}\dot{V} = \dot{V} - \xi N_\ast
 $$
 with a properly chosen $\xi \in \br$, whereas
$$
\|V\|\dot{V}_0 = \dot{V}_{V^\bot} = (\nabla_{\dot{\bldr}} V + \xi N^\ast)_{V^\bot} \overset{\text{(a)}}{=} (\nabla_{\dot{\bldr}} V )_{V^\bot} + \xi N^\ast,
$$
where (a) holds since $N_\ast$ is perpendicular to $V$.
Therefore,
\begin{gather*}
\|V\|\|W\|[N_\ast,\dot{V}_0,W_0] = [N_\ast,(\nabla_{\dot{\bldr}} V)_{V^\bot},W]
= \spr{(\nabla_{\dot{\bldr}} V)_{V^\bot}}{\mathscr{R}(-\pi/2)W}.
\end{gather*}
The proof is completed by handling $[N_\ast,V_0,\dot{W}_0]$ likewise. \epf
\begin{lemma}
\label{lem.merid}
Let $\bldb \in S$ and let $\vartheta\in (0,\pi/2]$ be the angle between the principal line $p_-(\bldb)$ and the line $l_{\bldb} = \{g \vec{\tau}(\bldb) + \bldb: g \in \br\}$ tangent to the meridian. The following inequalities hold:
\begin{equation}
\label{vartheta.ineq}
\Delta_\tau \leq 2L_N \quad \text{\rm and} \quad
\vartheta \geq \arcsin \sqrt{\frac{\Delta_{\tau}}{2L_N}},
\end{equation}
where $\Delta_\tau$ and $L_N$ are taken from \eqref{basic.ineq} and \eqref{lip.const}, respectively.
\end{lemma}
\pf
It suffices to note that
\begin{gather}
\label{ff.decom}
\ff{\bldb}{V}{V} = \varkappa_- \spr{V}{E_-}^2 + \varkappa_+ \spr{V}{E_+}^2,
\\
\nonumber
\varkappa_- + \Delta_\tau \overset{\text{\eqref{basic.ineq}}}{\leq} \varkappa_- \spr{\vec{\tau}}{E_-}^2 + \varkappa_+ \spr{\vec{\tau}}{E_+}^2
= \varkappa_- \cos^2 \vartheta + \varkappa_+ \sin^2 \vartheta
\\
\nonumber
\Rightarrow \Delta_\tau \leq (\varkappa_+ - \varkappa_-)\sin^2 \vartheta \overset{\text{\eqref{est.norm}}}{\leq} 2 L_N \sin^2 \vartheta \Rightarrow \text{\eqref{vartheta.ineq}. \epf}
\end{gather}
\begin{lemma}
\label{lem.liptau}
There exists $L_\tau >0$ such that
$\|\nabla_V\tau_\ast\| \leq L_\tau \|V\|$ for any vector $V \in \ts_{\bldb_\ast}[S(d_\ast)]$ and $\bldb_\ast \in S(d_\ast)$.
\end{lemma}
\pf
By \eqref{tau.ast}, $\vec{\tau}_\ast(\bldr)= \vec{\tau}_\lozenge[J_{-1}(\bldr)]$ for $\vec{\tau}_\lozenge(\bldb):= \frac{\vec{\tau}_\dagger(\bldb)}{\|\vec{\tau}_\dagger(\bldb)\|}$ and $\vec{\tau}_\dagger(\bldb):= J^\prime(\bldb) \vec{\tau}(\bldb)$. Hence
by Corollary~\ref{cor.smooth},
\begin{gather}
\label{est.n1}
\|\nabla_V \vec{\tau}_\ast\| \leq \|D_V  \vec{\tau}_\ast\| =\|D_{W}  \vec{\tau}_\lozenge\|,
\quad
\text{where}
\quad
W:= J_{-1}^\prime V \; \text{and} \;\|W\| \leq  \|V\|/\Delta_\varkappa;
\\
 1 = \|\vec{\tau}\| = \left\| J_{-1}^\prime \vec{\tau}_\dagger\right\| \leq \|\vec{\tau}_\dagger\|/\Delta_\varkappa,\quad \|D_{W}  \vec{\tau}_\lozenge\|
 \label{est.n2}
  = \frac{1}{\|\vec{\tau}_\dagger\|} \left\| D_W \vec{\tau}_\dagger - \spr{D_W \vec{\tau}_\dagger}{\vec{\tau}_\lozenge}\vec{\tau}_\lozenge\right\|
 \leq \frac{\|D_W \vec{\tau}_\dagger\|}{\Delta_\varkappa}.
 \end{gather}
 Now we invoke Assumption~\ref{ass.lipsch} and similarly see that $\|D_W \vec{b}\| \leq L_B \|W\|/\Delta_B^-$ for $\vec{b}:= \nabla B/\|\nabla B\|$. Since $\vec{\tau} = \pm \Big[E_- \spr{\vec{b}}{E_+} - E_+ \spr{\vec{b}}{E_-}\Big] $, we also get the following
\begin{gather}
\nonumber
\|D_W E_\pm\| \leq L_E \|W\|, \quad \|D_W \varkappa_\pm\| \leq L_\varkappa \|W\|,
\\
\nonumber
\Big\|D_WE_\pm \spr{\vec{b}}{E_\mp}\Big\| \leq \|D_W \vec{b}\| + \|D_W E_-\| + \|D_W E_+\|
\leq  (L_B/\Delta_B^- + 2 L_E)\|W\|,
\\
\nonumber
\|D_W \vec{\tau}\| \leq 2 (L_B/\Delta_B^- + 2 L_E)\|W\|,
\\
\nonumber
\|D_W(\varkappa_\pm \spr{E_\pm}{\vec{\tau}})\| \leq \{ L_\varkappa + |\varkappa_\pm| [L_E +  2 (L_B/\Delta_B^- + 2 L_E)]\} \|W\|
\\
\nonumber
\overset{\text{\eqref{est.norm}}}{\leq} [ L_\varkappa + L_N (5 L_E +  2 L_B/\Delta_B^-) ] \|W\|;
 \\
 \nonumber
 \vec{\tau}_\dagger := \vec{\tau} + d_\ast \mathscr{S}_{\bldb}\vec{\tau}
 =  \vec{\tau} + d_\ast \left[ \varkappa_- \spr{E_-}{\vec{\tau}}E_-+ \varkappa_+ \spr{E_+}{\vec{\tau}} E_+\right]
 \\
 \Rightarrow
\|D_W \vec{\tau}_\dagger\| \leq 2 (L_B/\Delta_B^- + 2 L_E)\|W\|
\label{est.n3}
+ 2 d_\ast [ L_\varkappa + L_N (6 L_E +  2 L_B/\Delta_B^-) ] \|W\|.
\end{gather}
The proof is completed by gathering \eqref{est.n1}---\eqref{est.n3}. \epf

\section{Proof of Proposition~\ref{prop.de}.}
\label{app.proofprop}
\setcounter{equation}{0}
Let $D_{\bcs}^\eta \subset \ts_{\bcs}(S)$ stands for the disc in the tangent plane $ \ts_{\bcs}(S)$ with a radius of $\eta>0$ centered at $\bcs$.
Formula \eqref{ff.decom} makes sense for any vector $V \in \br^3$. By using this, we extend the quadratic form $\ff{\bldb}{\cdot}{\cdot}$ and the associated symmetric bilinear form from the tangential plane to the entire space $\br^3$. If the tunnel is open, we put $S^0:= \{\bldb \in S: \bldb \,\text{is not in the $\delta_s$-vicinity of}\, \partial S\}$; otherwise, $S^0:=S$.
\begin{lemma}
\label{lemfform}
There exists $\eta>0$ such that for any point $\bcs \in S^0 $,
a patch $\mathscr{P}_{\bcs}(S)$ of $S$ around $\bcs$ is a graph of a $C^2$-smooth function $g_{\bcs}: D_{\bcs}^\eta \to \br$:
    \begin{equation}
    \label{patch}
    \mathscr{P}_{\bcs}(S) = \{\bldb = \Gamma_{\bcs}(\bldp):= \bldp + g_{\bcs}(\bldp) N(\bcs) : \bldp \in D_{\bcs}^\eta \}.
    \end{equation}
\end{lemma}
\pf
Let $\gamma$ be the normal section of $S$ by the plane $\mathpzc{P}$ that contains $\bcs$ and is coplanar with $N:= N(\bcs)$ and a unit vector $V \in \ts_{\bcs}(S)$.
There is a smooth function $f(\theta):\mathscr{E} \to \br, \mathscr{E}:=(-\ve_-,\ve_+), \ve_\pm >0$ such that its graph is a part of $\gamma$, i.e., $\gamma(\theta):=\bcs+\theta V + f(\theta) N \in \gamma, \theta \in \mathscr{E}$.
Since $f^\prime(0)=0$, reducing $\ve_\pm$, if necessary, ensures that for $L_N$ from \eqref{lip.const},
\begin{equation}
\label{def.dom}
L_N \|\gamma(\theta) - \gamma(0)\| + |f^\prime(\theta)|^2 < 1/2, \quad \|\gamma(\theta) - \gamma(0)\| < \delta_s/2 \qquad \forall \theta \in \mathscr{E}.
\end{equation}
From now on, we consider the maximal such an interval $\mathscr{E}$.
\par
The signed curvature $\varkappa(\theta)$ of $\gamma$ at point $\gamma(\theta)$ is
\begin{equation}
\label{Meusnier}
\varkappa(\theta) = \frac{f^{\prime\prime}(\theta)}{(1+|f^\prime(\theta)|^2)^{3/2}} \overset{\text{(a)}}{=}  \frac{\varkappa_{n}^S(\theta)}{\cos \beta}.
\end{equation}
Here $\beta$ is the angle between $V \times N$ and $V(\theta) \times N[\gamma(\theta)]$, whereas $\varkappa_{n}^S(\theta)$ is the normal curvature of $S$ at $\gamma(\theta)$ in direction of
\begin{equation}
\label{form.v}
V(\theta):=\gamma^\prime(\theta) = V+f^\prime(\theta)N,
\end{equation}
 and (a) holds by Meusnier's theorem \cite[p.~142]{Carmo76}. Here $\|V \times N\|=1$ and $\|V(\theta) \times N[\gamma(\theta)]\|=\|V(\theta)\| = \sqrt{1+|f^\prime(\theta)|^2}$. So $\|V(\theta)\|\cos \beta$ is the quadruple product of $V , N, V(\theta), N[\gamma(\theta)]$. By using formula (25) in \cite[p.~76]{GiWil01}
 $$
 \spr{A \times B}{C \times D} = \spr{A}{C} \spr{B}{D} - \spr{A}{D}\spr{B}{C},
 $$
 we see that
\begin{gather}
\sqrt{1+|f^\prime(\theta)|^2} \cos \beta
\nonumber
 = \spr{V}{V(\theta)}\spr{N}{N[\gamma(\theta)]} - \spr{V}{N[\gamma(\theta)]}\spr{V(\theta)}{N}
 \\
 \label{ugol.beta}
 \overset{\text{\eqref{form.v}}}{=} \spr{N}{N[\gamma(\theta)]} - f^\prime(\theta)\spr{V}{N[\gamma(\theta)]}
 \\
 \nonumber
 = 1 + \spr{N}{N[\gamma(\theta)]-N} - f^\prime(\theta)\spr{V-V(\theta)}{N[\gamma(\theta)]}
 \geq 1 - \|N[\gamma(\theta)]-N\| - |f^\prime(\theta)| \|V-V(\theta)\|
 \\
 \label{geq.half}
 \overset{\text{\eqref{lip.const},\eqref{form.v}}}{\geq} 1 - L_N \|\gamma(\theta) - \gamma(0)\| - |f^\prime(\theta)|^2 \overset{\text{\eqref{def.dom}}}{\geq} 1/2;
\\
\label{norm.curv}
\varkappa_{n}^S(\theta) = \frac{\ff{\gamma(\theta)}{\gamma^\prime(\theta)}{\gamma^\prime(\theta)}}{\|\gamma^\prime(\theta)\|^2}.
\end{gather}
Hence due to \eqref{Meusnier}, \eqref{form.v}, \eqref{geq.half}, and \eqref{norm.curv}
\begin{gather}
\frac{|f^{\prime\prime}(\theta)|}{(1+|f^\prime(\theta)|^2)^{3/2}} = |\varkappa(\theta)|
\leq 2 \frac{|\ff{\gamma(\theta)}{\gamma^\prime(\theta)}{\gamma^\prime(\theta)}|}{\|\gamma^\prime(\theta)\|^{3/2}}
\label{ddot.estimn}
\Rightarrow |f^{\prime\prime}(\theta)| \leq 2 |\ff{\gamma(\theta)}{\gamma^\prime(\theta)}{\gamma^\prime(\theta)}|.
\end{gather}
So thanks to \eqref{est.norm}, \eqref{ff.decom}, \eqref{ddot.estimn}, and Assumption~\ref{ass.lipsch},
\begin{gather}
\label{diff.ff}
|\ff{\bldb_2}{V}{V} - \ff{\bldb_1}{V}{V}|
\leq 2 (L_\varkappa + 2 L_E L_N) \|\bldb_2-\bldb_1\|\|V\|^2;
\\
\label{vikladka}
|f^{\prime\prime}(\theta)| \leq  2L_N \|\gamma^\prime(\theta)\|^2
=  2 L_N (1+|f^\prime(\theta)|^2).
\end{gather}
For some $\theta_+ =\theta_+(L_N)>0$,
the solution $q(\theta) = \tan (2 L_N \theta)$ of the Cauchy problem $q(0) =0$ for the ode
$
q^\prime = \lambda(q):= 2 L_N (1+q^2)
$
is defined on $[0,\theta_+]$ and
\begin{equation}
\label{q.ineq}
\int_0^{\theta_+} \sqrt{1+q(\varsigma)^2} \; d \varsigma \leq \delta_s/3, \quad
L_N \int_0^{\theta_+} \sqrt{1+q(\varsigma)^2} \; d \varsigma + q(\theta_+)^2 \leq 1/3.
\end{equation}
By \eqref{vikladka}, $q_\ast^\prime \leq \lambda(q_\ast)$ and $q_\ast(0)=0$ for $q_\ast(\theta):= \pm f^\prime(\theta), \pm f^\prime(-\theta)$. Hence by \cite[Th.~4.1, p.~26]{Hart82}, every of these $q_\ast(\theta)$'s does not exceed $q(\theta)$ whenever  $|\theta| < \min \{\ve_-,\ve_+\}$ and $|\theta| \leq \theta_+$, For such $\theta$, we thus have
\begin{equation}
\label{oc.g}
|f^\prime(\theta)| \leq q(|\theta|), \;
\|\gamma^\prime(\theta)\| \leq \sqrt{1+|q(|\theta|)|^2},
\;
\|\gamma^{\prime\prime}(\theta)\| \leq \lambda[q(|\theta|)],
\;
\|\gamma(\theta) - \gamma(0)\| \leq \int_0^{|\theta|} \sqrt{1+|q(|\varsigma|)|^2} \; d \varsigma.
\end{equation}
\par
We are going to show that $\min \{\ve_-,\ve_+\} > \theta_+$. Suppose to the contrary that $ \ve_\sigma=\min \{\ve_-,\ve_+\} \leq \theta_+$, where $\sigma = \pm$. Let $\sigma =+$, the case $\sigma=-$ is treated likewise.
For $\theta \in [0,\ve_+)$, we have
\begin{gather}
L_N\|\gamma(\theta) - \gamma(0)\| + |f^\prime(\theta)|^2
\overset{\text{\eqref{oc.g}}}{\leq} L_N \int_0^\theta \sqrt{1+q(\varsigma)^2}\;d \varsigma + q(\theta)^2
\label{fin.lim}
\overset{\text{\eqref{q.ineq}}}{\leq}  1/3 ,
\\
\label{partil.dal}
\|\gamma(\theta) - \gamma(0)\|
\overset{\text{\eqref{oc.g}}}{\leq} \int_0^\theta \sqrt{1+q(\varsigma)^2}\;d \varsigma \overset{\text{\eqref{q.ineq}}}{\leq}  \delta_s/3.
\end{gather}
Since $\|\gamma^\prime(\theta)\|$ and $\|\gamma^{\prime\prime}(\theta)\|$ stay bounded
as $\theta \to \ve_+-$ by \eqref{oc.g}, there exist $\bldb_{\lim} := \lim_{\theta \to \ve_+-}\gamma(\theta)$ and $W:= \lim_{\theta \to \ve_+-}\gamma^\prime(\theta)$ and the vector $W$ is not aligned with $N$. Meanwhile, $\bldb_{\lim} \not\in \partial S$ thanks to \eqref{partil.dal} since $\gamma(0) = \boldsymbol{c}$ is not in the $\delta_s$-vicinity of $\partial S$. It follows that the normal section $\gamma$ extends as a graph of a smooth function to the right of $\ve_+$. Meanwhile, \eqref{fin.lim} and \eqref{partil.dal} imply that
\eqref{def.dom} holds for $\theta=\ve_+$ and so
for $\theta > \ve_+, \theta \approx \ve_+$. These violate the definition of $\ve_+$ as an end of the maximal interval. The contradiction obtained proves that $\min \{\ve_-,\ve_+\} > \theta_+$.
\par
Now we put $\eta:=\theta_+, V_{\bldp} := \frac{\bldp - \bcs}{\|\bldp - \bcs\|}\; \forall \bldp \in D_{\bcs}^\eta, \bldp \neq \bcs$ and emphasize the dependence of $f(\cdot)$ on $V$ by adding the index $_V$ to $f$, thus obtaining $f_V(\cdot)$. Then the function
$$
g_{\bcs}(\bldp) :=
\begin{cases}
f_{V_{\bldp}}(\|\bldp - \bcs\|) & \text{whenever}\; \bldp \neq \bcs,
\\
0 & \text{if}\; \bldp = \bcs
\end{cases}
\qquad \forall \bldp \in D_{\bcs}^\eta
$$
is well-defined and meets \eqref{patch}. Since $\|N[\gamma(\theta)] - N\| \leq L_N \|\gamma(\theta) - \gamma(0)\| <1/2$ by \eqref{lip.const} and \eqref{def.dom}, the vectors $N$ and $N[\gamma(\theta)]$ are not aligned. It follows that the function $g_{\bcs}(\cdot)$ is $C^2$-smooth on $D_{\bcs}^\eta$. \epf
\begin{lemma}
\label{lem.grad}
After properly reducing $\eta>0$, if necessary, the following inequalities hold for any $\bldp \in D_{\bcs}^\eta, \bcs \in S^0$:
\begin{gather}
\label{ner.gradient}
\|\nabla g_{\bcs}(\bldp)\| \leq 2 L_N \|\bldp - \bcs\|,
\\
\label{ner.sama}
|g_{\bcs}(\bldp)| \leq  L_N \|\bldp - \bcs\|^2.
\end{gather}
\end{lemma}
\pf
For any unit vector $V \in \ts_{\bcs}(S)$ and $\bldp \in D_{\bcs}^\eta$, we have $\zeta(\theta) := \Gamma_{\bcs}[\bldp + \theta V] \in S \; \forall \theta \approx 0 \Rightarrow \ts_{\bldp}(S) \ni \zeta^\prime(0) = V + \spr{\nabla g_{\bcs}(\bldp)}{V}N(\bcs)$. So $\zeta^\prime(0)$ is normal to $N(\bldp)$ and
\begin{equation}
\label{start.from}
\spr{V}{N(\bldp)} = - \spr{N(\bcs)}{N(\bldp)} \spr{\nabla g_{\bcs}(\bldp)}{V}.
\end{equation}
Here $|\spr{V}{N(\bldp)}| = |\spr{V}{N(\bldp)-N(\bcs)}| \leq \|N(\bldp)-N(\bcs)\|\leq L_N \|\bldp - \bcs\|$ by \eqref{lip.const}. Meanwhile, $\spr{N(\bcs)}{N(\bldp)} = 1+\spr{N(\bcs)}{N(\bldp)-N(\bcs)}$, where $|\spr{N(\bcs)}{N(\bldp)-N(\bcs)}| \leq \|N(\bldp)-N(\bcs)\| \leq L_N \|\bldp - \bcs\|$. So by properly reducing $\eta>0$, if necessary, we can ensure that
$\spr{N(\bcs)}{N(\bldp)} \geq 1/2$. Then \eqref{start.from} implies that $|\spr{\nabla g_{\bcs}(\bldp)}{V}| \leq 2 L_N \|\bldp - \bcs\|$. Maximization over $V \in \ts_{\bcs}, \|V\|=1$ yields \eqref{ner.gradient}.
Also,
\begin{gather*}
|g_{\bcs}(\bldp)| = \left| \int_0^1 \spr{\nabla g_{\bcs}[\theta\bldp + (1-\theta)\bcs]}{\bldp - \bcs} \; d \theta \right|
\overset{\text{\eqref{ner.gradient}}}{\leq} 2 L_N \|\bldp - \bcs\|^2 \int_0^1 \theta \; d \theta \Rightarrow \text{\eqref{ner.sama} . \epf}
\end{gather*}
\par
Based on \eqref{patch} and \eqref{ner.sama}, we see that
\begin{equation}
\label{sama1}
\|\Gamma_{\bcs}(\bldp) - \bcs\| \leq \|\bldp - \bcs\| + L_N \|\bldp - \bcs\|^2.
\end{equation}
\begin{lemma}
\label{lem.secondder}
There is a non-decaying function $\zeta(\cdot)$ such that $\zeta(\varrho) \to 0$ as $\varrho \to 0+$ and
the Hessian $g^{\prime\prime}_{\bcs}$ obeys the estimate
$$
\|g^{\prime\prime}_{\bcs}(\bldp) - g^{\prime\prime}_{\bcs}(\bcs)\| \leq \zeta(\|\bldp - \bcs\|) \qquad \forall \bldp \in D_{\bcs}^\eta, \bcs \in S^0.
$$
\end{lemma}
\pf
In $\br^3$, we pick a Cartesian coordinate system centered at $\bcs$ so that its $x$- and $y$-axes lie in the tangent plane $\ts_{\bcs}(S)$ and the $z$-axis is co-directed with $N(\bcs)$.
We also identify any point $\bldp \in \ts_{\bcs}(S)$ with the pair $(x,y)$ of its coordinates.
Since the map $\Gamma:=\Gamma_{\bcs} =[x,y,z(x,y)], z(x,y) := g_{\bcs}(x,y)$ is a coordinate chart on $S$, we have \cite[p.~154]{Carmo76}
\begin{gather}
\label{dopil}
\ff{\Gamma(\bldp)}{V_{\bldp}}{V_{\bldp}} = A(\bldp) (dx)^2 + 2 B(\bldp) dx dy + C(\bldp) (dy)^2,
\\
\nonumber
\text{where}\quad V_{\bldp}:=  \Gamma_x^\prime(\bldp) dx + \Gamma_y^\prime(\bldp) dy = [dx, dy, z^\prime_x dx + z^\prime_y dy],
\\
\nonumber
A(\bldp)= \spr{N_{\Gamma(\bldp)}}{\Gamma^{\prime\prime}_{xx}(\bldp)} = \spr{N_{\Gamma(\bldp)}}{N_{\bcs}} z^{\prime\prime}_{xx}(\bldp),
\\
\nonumber
B(\bldp)= \spr{N_{\Gamma(\bldp)}}{\Gamma^{\prime\prime}_{xy}(\bldp)}= \spr{N_{\Gamma(\bldp)}}{N_{\bcs}} z^{\prime\prime}_{xy}(\bldp),
\\
\nonumber
C(\bldp)= \spr{N_{\Gamma(\bldp)}}{\Gamma^{\prime\prime}_{yy}(\bldp)} = \spr{N_{\Gamma(\bldp)}}{N_{\bcs}} z^{\prime\prime}_{yy}(\bldp).
\end{gather}
After reducing $\eta$ so that $\eta <1$, we have for $\bldp \in D_{\bcs}^\eta$,
\begin{gather*}
|\ff{\Gamma(\bldp)}{V_{\bldp}}{V_{\bldp}} - \ff{\Gamma(\bcs)}{V_{\bcs}}{V_{\bcs}}|
\leq |\ff{\Gamma(\bldp)}{V_{\bldp}}{V_{\bldp}} - \ff{\Gamma(\bcs)}{V_{\bldp}}{V_{\bldp}}|
+
|\ff{\Gamma(\bcs)}{V_{\bldp}}{V_{\bldp}} - \ff{\Gamma(\bcs)}{V_{\bcs}}{V_{\bcs}}|
\\
\overset{\text{\eqref{diff.ff}}}{\leq} 2 (L_\varkappa + 2 L_E L_N) \|\Gamma(\bldp)-\Gamma(\bcs)\|\|V_{\bldp}\|^2
+ |\ff{\Gamma(\bcs)}{V_{\bldp}- V_{\bcs}}{V_{\bldp}- V_{\bcs}}| + 2 |\ff{\Gamma(\bcs)}{V_{\bldp}- V_{\bcs}}{V_{\bcs}}|
\\
\overset{\text{\eqref{est.norm},\eqref{ff.decom},\eqref{sama1}}}{\leq}
 2 (L_\varkappa + 2 L_E L_N)\eta (1+ L_N \eta) \|V_{\bldp}\|^2
 + L_N [\|V_{\bldp} - V_{\bcs}\|^2 + 2\|V_{\bldp} - V_{\bcs}\| \|V_{\bcs}\|].
\end{gather*}
By noting that $\nabla g_{\bcs}(\bcs)=0$ due to \eqref{ner.gradient}, we have $\|V_{\bldp} - V_{\bcs}\| = |[z^\prime_x(\bldp) - z^\prime_x(\bcs)]dx + [z^\prime_y(\bldp) - z^\prime_y(\bcs)] dy|
\leq || \nabla g_{\bcs}(\bldp)\| \xi_d$, where $ \xi_d:= \sqrt{(dx)^2+(dy)^2}$. So by Lemma~\ref{lem.grad}, $\|V_{\bldp} - V_{\bcs}\| \leq 2 L_N \eta \xi_d$. Meanwhile $\|V_{\bcs}\|=\xi_d$. Thus
\begin{gather*}
|\ff{\Gamma(\bldp)}{V_{\bldp}}{V_{\bldp}} - \ff{\Gamma(\bcs)}{V_{\bcs}}{V_{\bcs}}| \leq M_{\mathbf{II}} \eta \,\xi_d^2, \qquad \text{where}
\\
M_{\mathbf{II}}:=  2 (L_\varkappa + 2 L_E L_N)(1 +L_N) [1+2 L_N  ]^2
 + 4 L_N^2 (L_N  + 1 )  .
\end{gather*}
 Here $\ff{\Gamma(\bldp)}{V_{\bldp}}{V_{\bldp}} - \ff{\Gamma(\bcs)}{V_{\bcs}}{V_{\bcs}}= [A(\bldp) - A(\bcs)](dx)^2 + 2 [B(\bldp) -B(\bcs)] dx dy + [C(\bldp) -C(\bcs)] (dy)^2$. Hence
 \begin{gather*}
 |A(\bldp) - A(\bcs)| \leq M_{\mathbf{II}} \eta, \quad |C(\bldp) - C(\bcs)| \leq M_{\mathbf{II}} \eta,
  \quad
  |B(\bldp) - B(\bcs)| \leq  M_{\mathbf{II}} \eta.
 \end{gather*}
At the same time,
\begin{gather*}
A(\bldp) - A(\bcs) = \spr{N_{\Gamma(\bldp)}}{N_{\bcs}} z^{\prime\prime}_{xx}(\bldp) -  z^{\prime\prime}_{xx}(\bcs)
=
 \spr{N_{\Gamma(\bldp)}}{N_{\bcs}} [z^{\prime\prime}_{xx}(\bldp) -  z^{\prime\prime}_{xx}(\bcs)]
 + [\spr{N_{\Gamma(\bldp)}}{N_{\bcs}}-1] z^{\prime\prime}_{xx}(\bcs);
 \\
 |\spr{N_{\Gamma(\bldp)}}{N_{\bcs}}-1| = |\spr{N_{\Gamma(\bldp)}-N_{\bcs}}{N_{\bcs}}|
 \overset{\text{\eqref{lip.const}}}{\leq} L_N \|\Gamma(\bldp) - \bcs\| \overset{\text{\eqref{sama1}}}{\leq}
 L_N [\|\bldp - \bcs\| + L_N \|\bldp - \bcs\|^2]
\\
 \leq L_N [1 + L_N \eta]\eta, \quad  \spr{N_{\Gamma(\bldp)}}{N_{\bcs}}  \geq 1 - L_N \eta [1 + L_N \eta] \geq 1/2
\end{gather*}
if $\eta < (\sqrt{3}-1)/(2L_N)$. Meanwhile \eqref{dopil} implies that
\begin{equation}
\label{rav.vtor}
\spr{g^{\prime\prime}_{\bcs}(\bcs)V}{V}=  \ff{\bcs}{V}{V}
\end{equation}
and so $|\spr{z^{\prime\prime}(\bcs)V}{V}|=|\spr{g^{\prime\prime}_{\bcs}(\bcs)V}{V}| = |\ff{\bcs}{V}{V}| = |\spr{\mathscr{S}_{\bcs}V}{V}| \leq L_N\|V\|^2$ by \eqref{est.norm}. Hence $|z^{\prime\prime}_{xx}(\bcs)| \leq L_N,
|z^{\prime\prime}_{yy}(\bcs)| \leq L_N, |z^{\prime\prime}_{xy}(\bcs)| \leq L_N$. As a result,
\begin{gather*}
|z^{\prime\prime}_{xx}(\bldp) -  z^{\prime\prime}_{xx}(\bcs)|
= \frac{|A(\bldp)- A(\bcs) - [\spr{N_{\Gamma(\bldp)}}{N_{\bcs}}-1] z^{\prime\prime}_{xx}(\bcs)|}{\spr{N_{\Gamma(\bldp)}}{N_{\bcs}}}
\leq 2 (M_{\mathbf{II}} + L_N^2 +L_N^3 \eta)\eta.
\end{gather*}
Similarly $|z^{\prime\prime}_{yy}(\bldp) -  z^{\prime\prime}_{yy}(\bcs)| \leq 2 (M_{\mathbf{II}} + L_N^2+L_N^3 \eta)\eta, |z^{\prime\prime}_{xy}(\bldp) -  z^{\prime\prime}_{xy}(\bcs)| \leq 2 (M_{\mathbf{II}} + L_N^2+L_N^3 \eta)\eta$.
It remains to note that for a given $\bldp$ from the basic disk $D_{\bcs}^\eta$, the above estimates remain true with $\eta$ being artificially reduced and made arbitrarily close to $\|\bldp - \bcs\|$. \epf
\begin{corr}
\label{cor.basicin}
The following relations hold for any $\bcs \in S^0$, $\bldp \in D_{\bcs}^\eta$:
\begin{gather*}
g_{\bcs}(\bldp) = 1/2 \ff{\bcs}{\bldp-\bcs}{\bldp-\bcs} + \omega_{\bcs}(\bldp), \quad \text{\rm where}
\\
|\omega_{\bcs}(\bldp)| \leq \zeta(\|\bldp - \bcs\|) \|\bldp - \bcs\|^2/2, \\ \|\nabla \omega_{\bcs}(\bldp)\|
\leq \zeta(\|\bldp - \bcs\|) \|\bldp - \bcs\|,
\\
\|\omega_{\bcs}^{\prime\prime}(\bldp)\|
\leq \zeta(\|\bldp - \bcs\|).
\end{gather*}
\end{corr}
The last inequality is immediate from Lemma~\ref{lem.secondder} and \eqref{rav.vtor}. Then putting $\bldp(\theta):= (\bldp-\bcs)\theta + \bcs$, we have
\begin{gather*}
\left\| \nabla \omega_{\bcs}(\bldp) \right\| = \left\| \int_0^1 \omega_{\bcs}^{\prime\prime}[\bldp(\theta)][\bldp - \bcs]\; d \theta \right\|
\leq \|\bldp - \bcs\|\int_0^1 \|\omega_{\bcs}^{\prime\prime}[\bldp(\theta)]\| \; d \theta \leq \zeta(\|\bldp - \bcs\|)\|\bldp - \bcs\|;
\\
\left|\omega_{\bcs}(\bldp) \right| = \left| \int_0^1 \spr{\nabla \omega_{\bcs}[\bldp(\theta)]}{\bldp - \bcs}\; d \theta \right|
\leq \|\bldp - \bcs\| \int_0^1 \|\nabla \omega_{\bcs}[\bldp(\theta)]\|\; d \theta \leq \frac{\zeta(\|\bldp - \bcs\|)}{2}\|\bldp - \bcs\|^2 .
\end{gather*}
\par
\begin{figure}
\centering
\scalebox{0.3}{\includegraphics{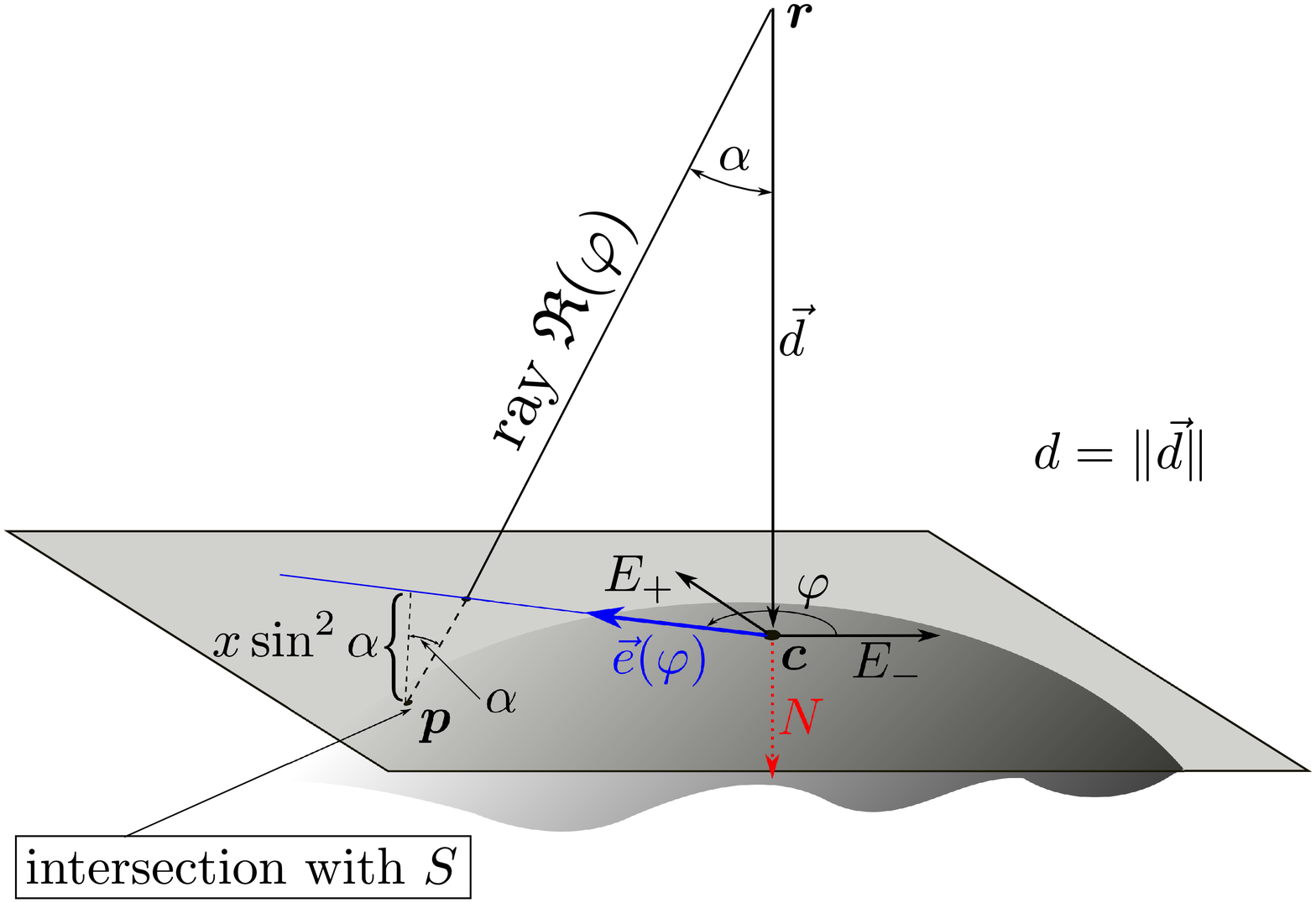}}
\caption{}
\label{definx.fig}
\end{figure}
Now we consider a point $\bldr \in Z_{\text{op}}^a$ and its projection $\bcs = \bcs(\bldr):= \bpi(\bldr)$ onto $S$; then $\bcs \in S^0$ by the definition of $Z_{\text{op}}^a$. The rays $R$ emitted from $\bldr$ at a given angle $\alpha \in [0, \alpha_s]$ to $\vec{d} = \bcs - \bldr$ are parametrized $\mathfrak{R}=\mathfrak{R}(\varphi)$ by angle $\varphi$ so that $\mathfrak{R}(\varphi)$ passes through the point $\bcs+ d \vec{e}(\varphi) \tan \alpha $, where $d = \|\vec{d}\|$ is the distance from $\bldr$ to $S$ and $\vec{e}(\varphi) := E_-(\bcs) \cos \varphi +  E_+(\bcs)\sin \varphi$; see Fig.~\ref{definx.fig}.
The distance $d(\alpha,\varphi)$ from $\bldr$ to the point $\bldp$ of intersection of $\mathfrak{R}(\varphi)$ with the surface $S$
is given by $d(\alpha,\varphi) = (d+x \sin^2\alpha )/\cos \alpha$, where $x=x(\alpha,\varphi)$ is defined in Fig.~\ref{definx.fig} and scaling of the distance from $\bldp$ to the tangent plane by the factor $\sin^2 \alpha$ is introduced for further convenience. It follows that maximization of $d(\alpha,\varphi)$ over $\varphi$ is equivalent to maximization of $x$. Meanwhile, $x$ is the root
of the equation
\begin{equation}
\label{definition.x}
x = \lambda(x,\varphi|\alpha) := - \frac{g_{\bcs}[\bcs+ (d+x\sin^2\alpha) \vec{e}(\varphi) \tan \alpha]}{\sin^2\alpha}.
\end{equation}
For this equation to be well-posed, the argument of $g_{\bcs}(\cdot)$ should be in the domain $D_{\bcs}^\eta$ of definition of $g_{\bcs}(\cdot)$, which is equivalent to
\begin{equation}
\label{x.interval}
x \in I:=\left[ - \eta \frac{\cos\alpha}{\sin^3\alpha} - \frac{d}{\sin^2\alpha}; \eta \frac{\cos\alpha}{\sin^3\alpha} - \frac{d}{\sin^2\alpha}\right].
\end{equation}
\begin{lemma}
\label{lemma.smooths}
There is $\alpha_1 = \alpha_1(\eta) \in (0,\alpha_s]$ such that whenever $\alpha \leq \alpha_1$, equation \eqref{definition.x} has an unique root on the interval \eqref{x.interval} and this root smoothly depends on $\varphi$ and $\bldr \in Z_{\text{op}}^a$.
\end{lemma}
\pf
We first note that for $x \in I$,
\begin{gather*}
 |\lambda^\prime_x| \cos \alpha = |\spr{\nabla g_{\bcs}[\bcs+ (d+x\sin^2\alpha) \vec{e}(\varphi) \tan \alpha]}{ \vec{e}(\varphi)}\sin\alpha|
 \\
 \overset{\text{\eqref{ner.gradient}}}{\leq} 2L_N  |d+x\sin^2\alpha| \tan \alpha \sin\alpha
 \overset{\text{\eqref{x.interval}}}{\leq} 2L_N \eta \sin\alpha.
\end{gather*}
So $|\lambda^\prime_x| \leq 1/2\; \forall x \in I, \alpha \in (0,\alpha_1]$ if $\alpha_1$ is chosen small enough. This implies that the root $x \in I$ of \eqref{definition.x} (if exists) is unique. At the right end of $I$, the r.h.s. of \eqref{definition.x} equals $-g_{\bcs}[\bcs+\eta \vec{e}(\varphi)]\sin^{-2}\alpha$ and is bounded over $\varphi$ due to \eqref{ner.sama}, whereas the l.h.s. is
$[\eta \frac{\cos\alpha}{\sin\alpha} - d] \sin^{-2}\alpha \geq [\eta \frac{\cos\alpha}{\sin\alpha} - d_+] \sin^{-2}\alpha$, where $\eta \frac{\cos\alpha}{\sin\alpha} - d_+ \to \infty$ as $\alpha \to 0+$.
So by reducing $\alpha_1>0$, if necessary, it can be ensured that r.h.s. $<$  l.h.s. at the right end of $I$ for all $\varphi$ and $\bldr \in Z_{\text{op}}^a$. Similarly, it can be ensured that r.h.s. $>$ l.h.s. at the left end. Hence the root $x \in I$ of equation \eqref{definition.x} does exist. By the implicit function theorem (applied to the function $x-\lambda(x,\varphi|\alpha)$) and Lemma~\ref{lem.cuvb}, this root smoothly depends on $\varphi$ and $\bldr \in Z_{\text{op}}^a$. \epf
\begin{lemma}
\label{lem.unifcon}
As $\alpha \to 0+$, the function $\lambda(x,\varphi|\alpha)$ and its first and second derivatives in $x$ and $\varphi$ converge to
the function $\lambda(x,\varphi|0)= - d^2 \ff{\bcs}{\vec{e}(\varphi)}{\vec{e}(\varphi)} /2$ and its respective derivatives uniformly over $\bldr \in Z_{\text{op}}^a, \varphi \in \br$, and $x \in [-\ov{x}, \ov{x}]$ for any $\ov{x}>0$.
\end{lemma}
\pf
For $\bldp := \bcs+ (d+x\sin^2\alpha) \vec{e}(\varphi) \tan \alpha, x \in [-\ov{x}, \ov{x}]$, and $k:= d_+ +\ov{x}$, we have $\|\bldp - \bcs\|=|x \sin^2\alpha +d| \tan \alpha \leq k \tan \alpha$.
Thanks to Corollary~\ref{cor.basicin},
\begin{gather*}
\lambda(x,\varphi|\alpha) = A+B, \quad \text{where}
\\
A:=-\frac{\ff{\bcs}{(d+x\sin^2\alpha) \vec{e}(\varphi) }{(d+x\sin^2\alpha) \vec{e}(\varphi) }}{2\cos^2\alpha},
\\
B:= -\frac{\omega_{\bcs}(\bldp)}{\sin^2\alpha};
 \quad
|B| \leq  \zeta ( k \tan \alpha) \frac{k^2 }{2 \cos^2\alpha},
\\
|B^\prime_x| \leq \left|\spr{\nabla \omega_{\bcs}(\bldp)}{ \vec{e}(\varphi)} \right| \tan \alpha
\leq k \zeta ( k \tan \alpha)  \tan^2 \alpha,
\\
|B^\prime_\varphi| \leq \frac{|\spr{\nabla \omega_{\bcs}(\bldp)}{ \vec{e}^{\,\prime}(\varphi)}||d+x\sin^2\alpha|}{\sin\alpha\cos \alpha}
\leq k \zeta ( k \tan \alpha)/ \cos^2 \alpha,
\\
|B^{\prime\prime}_{xx}| \leq |\spr{\omega_{\bcs}^{\prime\prime}(\bldp)\vec{e}(\varphi)}{\vec{e}(\varphi)}| \tan^2 \alpha \sin^2\alpha
\leq \zeta ( k \tan \alpha)  \tan^2 \alpha \sin^2\alpha ,
\\
|B^{\prime\prime}_{x\varphi}| \leq |\spr{\nabla \omega_{\bcs}(\bldp)}{\vec{e}^{\,\prime}(\varphi)}|\tan \alpha
+ |\spr{\omega_{\bcs}^{\prime\prime}(\bldp)\vec{e}(\varphi)}{\vec{e}^{\,\prime}(\varphi)}|  |d+x \sin^2\alpha|  \tan^2\alpha
\leq
2 k \zeta ( k \tan \alpha) \tan^2 \alpha ,
\\
|B^{\prime\prime}_{\varphi\varphi}| \leq \frac{|\spr{\nabla \omega_{\bcs}(\bldp)}{\vec{e}^{\,\prime\prime}(\varphi)}||d+x\sin^2\alpha|}{\sin \alpha \cos \alpha}
+ \frac{|\spr{\omega_{\bcs}^{\prime\prime}(\bldp) \vec{e}^{\, \prime}(\varphi)}{\vec{e}^{\, \prime}(\varphi)}|(d+x\sin^2\alpha)^2 }{\cos^2\alpha}
\leq 2 \frac{k^2 \zeta ( k \tan \alpha)}{\cos^2 \alpha}.
 \end{gather*}
Thus $B$ and its first and second derivatives uniformly converge to $0$ as $\alpha \to 0+$. It remains to note that $A$ and its first and second derivatives uniformly converge to the function
$\lambda(x,\varphi|0)$ and its respective derivatives since $\ff{\bcs}{\cdot}{\cdot}$ is a quadratic form. \epf
\par
In the light of Lemmas~\ref{lemma.smooths} and \ref{lem.unifcon}, the implicit function theorem applied to the function $x - \lambda(x,\varphi|\alpha)$ yields the following properties of the root $x=x(\varphi)$ of equation \eqref{definition.x} (which also depends on $\bldr$ through $d$ and $\bcs$).
\begin{corr}
\label{corr.final}
The function $x(\varphi)$ is $C^2$-smooth. This function and its first and second derivatives are continuous in $\bldr$ and converge to the root $y(\varphi) = - d^2 \ff{\bcs}{\vec{e}(\varphi)}{\vec{e}(\varphi)} /2$ of the equation $y = \lambda(y,\varphi|0)$ and its respective derivatives uniformly over $\varphi \in \br, \bldr \in Z_{\text{op}}^a$ as $\alpha \to 0+$.
\end{corr}
\par
{\bf PROOF OF PROPOSITION~\ref{prop.de}.}
Due to the definition of $\vec{e}(\varphi)$ (see Fig.~\ref{definx.fig}) and \eqref{ff.decom},
\begin{gather*}
4 y(\varphi) = - 2 d^2 [\varkappa_-(\bcs) \cos^2\varphi + \varkappa_+(\bcs) \sin^2\varphi]
= - d^2 [\varkappa_-(\bcs)+\varkappa_+(\bcs)]   + d^2 [\varkappa_+(\bcs) - \varkappa_-(\bcs)]\cos(2 \varphi) .
\end{gather*}
Since $\Delta_\tau \leq \varkappa_+ - \varkappa_- \leq 2L_N$ due to \eqref{basic.ineq} and \eqref{est.norm}, we have
\begin{gather*}
y^\prime(\varphi) = 0 \Leftrightarrow \varphi =0,  \pm \pi/2, \pi \mod  2 \pi,
\\
y^{\prime\prime}(\varphi) \leq - d_-^2 \Delta_\tau <0 \;\text{for}\; \varphi = 0, \pi  \mod 2 \pi,
\\
y^{\prime\prime}(\varphi) \geq  d_-^2 \Delta_\tau > 0 \;\text{for}\; \varphi =  \pm \pi/2  \mod  2 \pi ,
\\
|y^{\prime\prime}(\varphi)| \leq  d^2 |\varkappa_+-\varkappa_-| \leq 2 d_+^2 L_N,
\quad
|y^{\prime\prime\prime}(\varphi)| \leq 4 d_+^2 L_N.
\end{gather*}
So putting $\delta:= \frac{d_-^2 \Delta_\tau}{8d_+^2 L_N}$, we have $y(0)=y(\pi)$,
\begin{gather*}
y^{\prime\prime}(\varphi) \leq - d_-^2 \Delta_\tau/2 <0 \;
\forall \varphi \in (-\delta, \delta)\cup (\pi-\delta,\pi+\delta)  \mod 2 \pi,
\\
y^{\prime\prime}(\varphi) \geq  d_-^2 \Delta_\tau/2 > 0 \; \forall \varphi \in (\pm \pi/2 - \delta, \pm \pi/2+\delta)  \mod 2 \pi,
\\
|y^\prime(\varphi)| \geq d_-^2 \delta \Delta_\tau/2 \quad \text{whenever $\varphi$ is outside all of the listed intervals}.
\end{gather*}
Corollary~\ref{corr.final} ensures that the following relation hold for all small enough $\alpha$:
\begin{gather*}
x^{\prime\prime}(\varphi) \leq - d_-^2 \Delta_\tau/4 <0 \;
\forall \varphi \in (-\delta, \delta)\cup (\pi-\delta,\pi+\delta)  \mod 2 \pi,
\\
x^{\prime\prime}(\varphi) \geq  d_-^2 \Delta_\tau/4 > 0 \; \forall \varphi \in (\pm \pi/2 - \delta, \pm \pi/2+\delta)  \mod 2 \pi,
\\
\text{and}\quad |x^\prime(\varphi)| \geq d_-^2 \delta \Delta_\tau/4 \quad \text{whenever $\varphi$ is outside all of the listed intervals}.
\end{gather*}
By applying the implicit function theorem to the equation $x^\prime(\varphi)=0$ and reducing $\delta$, if necessary, we see that
this equation has a single root on every of the following $(\mathrm{mod}\, 2\pi)$-intervals $(-\delta, \delta), (\pi-\delta,\pi+\delta), (- \pi/2 - \delta, - \pi/2+\delta), (\pi/2 - \delta, \pi/2+\delta)$, this root continuously depends on $\bldr \in Z_{\text{op}}^a$ and uniformly goes to the center of the respective interval as $\alpha \to 0+$. Meanwhile, the last displayed formulas assure that there are no other roots and $x(\varphi)$ attains its local maxima at two points: at the roots $\varphi_0 \in (-\delta,\delta)$ and $\varphi_\pi \in (\pi -\delta,\pi+\delta)$. So MDPBE is well posed. It remains to note that it returns the line that goes in the direction of $\varphi_\ast:=1/2(\varphi_0 + \varphi_\pi - \pi)$, where $\varphi_\ast \to 0$ as $\alpha \to 0+$ uniformly over $\bldr \in Z_{\text{op}}^a$ and $\varphi=0$ corresponds to the principal direction $p_-(\bcs)$. \epf

\end{document}